\documentstyle[txmac,a4,%leqno,%
amssymb,%
case,%
twoside,%
nocaphead,%
epsf,%
varthm,%
rotate,%
myrot,%
mypic,times,mathptm]{article}

\advance\oddsidemargin by -0.9cm
\advance\evensidemargin by -0.9cm
\advance\textwidth by 1.8cm
\def\mynewtheo#1#2{%
\newtheorem{@#1}{#2}%
\newenvironment{#1}{\begin{@#1}\rm}{\end{@#1}}}

\mynewtheo{lemma}{Lemma}
\mynewtheo{exer}{Exercise}
\mynewtheo{theo}{Theorem}
\mynewtheo{rem}{Remark}
\mynewtheo{defi}{Definition}
\mynewtheo{conj}{Conjecture}
\mynewtheo{corr}{Corollary}
\mynewtheo{prop}{Proposition}
\mynewtheo{question}{Question}
\mynewtheo{exam}{Example}

\newenvironment{theorem}{\begin{theo}}{\end{theo}}
\newenvironment{conjecture}{\begin{conj}}{\end{conj}}

\begin{document}

\makeatletter

% typing `equation' breaks my fingers!
\newenvironment{eqn}{\begin{equation}}{\end{equation}\@ignoretrue}

% eqlabel=(#1)
\newenvironment{myeqn*}[1]{\begingroup\def\@eqnnum{\reset@font\rm#1}%
\xdef\@tempk{\arabic{equation}}\begin{equation}\edef\@currentlabel{#1}}
{\end{equation}\endgroup\setcounter{equation}{\@tempk}\ignorespaces}

% eqlabel=#1
\newenvironment{myeqn}[1]{\begingroup\let\eq@num\@eqnnum
\def\@eqnnum{\bgroup\let\r@fn\normalcolor % an extremely UGLY hack !!!
\def\normalcolor####1(####2){\r@fn####1#1}%
%\show\reset@font
\eq@num\egroup}%
\xdef\@tempk{\arabic{equation}}\begin{equation}\edef\@currentlabel{#1}}
{\end{equation}\endgroup\setcounter{equation}{\@tempk}\ignorespaces}

% eqlabel=(eqnnr) \qed
\newenvironment{myeqn**}{\begin{myeqn}{(\arabic{equation})\es\es\mbox{\qed}}\edef\@currentlabel{\arabic{equation}}}
{\end{myeqn}\stepcounter{equation}}

\newcommand{\mybin}[2]{\text{$\Bigl(\begin{array}{@{}c@{}}#1\\#2%
\end{array}\Bigr)$}}
\newcommand{\mybinn}[2]{\text{$\biggl(\begin{array}{@{}c@{}}%
#1\\#2\end{array}\biggr)$}}

\def\overtwo#1{\mbox{\small$\mybin{#1}{2}$}}
\newcommand{\mybr}[2]{\text{$\Bigl\lfloor\mbox{%
\small$\displaystyle\frac{#1}{#2}$}\Bigr\rfloor$}}
\def\mybrtwo#1{\mbox{\mybr{#1}{2}}}

\def\myfrac#1#2{\raisebox{0.2em}{\small$#1$}\!/\!\raisebox{-0.2em}{\small$#2$}}

\def\myeqnlabel{\bgroup\@ifnextchar[{\@maketheeq}{\immediate
\stepcounter{equation}\@myeqnlabel}}
% label someth outside of eqns
% NOTE: this label does NOT work to be referenced by \pageref !

\def\@maketheeq[#1]{\def\theequation{#1}\@myeqnlabel}

\def\@myeqnlabel#1{%
%\immediate\stepcounter{equation}%
{\edef\@currentlabel{\theequation}
\label{#1}\enspace\eqref{#1}}\egroup}

\def\rato#1{\hbox to #1{\rightarrowfill}}
\def\arrowname#1{{\enspace
\setbox7=\hbox{F}\setbox6=\hbox{%
\setbox0=\hbox{\footnotesize $#1$}\setbox1=\hbox{$\to$}%
\dimen@\wd0\advance\dimen@ by 0.66\wd1\relax
$\stackrel{\rato{\dimen@}}{\copy0}$}%
\ifdim\ht6>\ht7\dimen@\ht7\advance\dimen@ by -\ht6\else
\dimen@\z@\fi\raise\dimen@\box6\enspace}}

\def\namedarrow#1{{\enspace
\setbox7=\hbox{F}\setbox6=\hbox{%
\setbox0=\hbox{\footnotesize $#1$}\setbox1=\hbox{$\to$}%
\dimen@\wd0\advance\dimen@ by 0.66\wd1\relax
$\stackrel{\rato{\dimen@}}{\copy0}$}%
\ifdim\ht6>\ht7\dimen@\ht7\advance\dimen@ by -\ht6\else
\dimen@\z@\fi\raise\dimen@\box6\enspace}}
   
\def\epsfs#1#2{{\epsfxsize#1\relax\epsffile{#2.eps}}}

\author{A. Stoimenow\footnotemark[1]\\[2mm]
\small Department of Mathematics,\\
\small University of Toronto,\\
\small Canada M5S 3G3,\\
\small e-mail: {\tt stoimeno@math.toronto.edu},\\
\small URL: {\hbox{\tt http://www.math.toronto.edu/stoimeno/}}
}

{\def\thefootnote{\fnsymbol{footnote}}
\footnotetext[1]{Supported by a DFG postdoc grant.}
}

\title{\large\bf \uppercase{A property of the skein polynomial with}\\[2mm]
\uppercase{an application to contact geometry}\\[4mm]
%\phantom
{\small\it This is a preprint. I would be grateful
for any comments and corrections!}
}

\date{\large Current version: \today\ \ \ First version:
\makedate{19}{6}{2000}}

\maketitle

\makeatletter

\let\point\pt
\let\ay\asymp
\let\pa\partial
\let\al\alpha
\let\be\beta
\let\Gm\Gamma
\let\gm\gamma
\let\de\delta
\let\dl\delta
\let\eps\epsilon
\let\lm\lambda
\let\Lm\Lambda
\let\sg\sigma
\let\vp\varphi
\let\om\omega
\let\fa\forall

\def\ssim{\stackrel{\ds \sim}{\vbox{\vskip-0.2em\hbox{$\scriptstyle *$}}
}}

\let\sm\setminus
\let\tl\tilde
\def\ncap{\not\mathrel{\cap}}
\def\pr{\text{\rm pr}\,}
\def\lra{\longrightarrow}
\def\Lra{\leftrightarrow}
\def\llra{\longleftrightarrow}
\def\so{\Rightarrow}
\def\So{\Longrightarrow}
\let\ds\displaystyle
\def\bt{\bar t_2}
\def\Md{\max\deg}
\def\md{\min\deg}
\def\spn{{\mathop {\operator@font span}}}
\def\mc{\min\cf}
\def\cf{\mathop{\operator@font cf}}

\let\reference\ref

\long\def\@makecaption#1#2{%
   % \tm
   \vskip 10pt
   {\let\label\@gobble
   \let\ignorespaces\@empty
   \xdef\@tempt{#2}%
   %\typeout{`#2'}%
   }%
   \ea\@ifempty\ea{\@tempt}{%
   \setbox\@tempboxa\hbox{%
      \fignr#1#2}%
      }{%
   \setbox\@tempboxa\hbox{%
      {\fignr#1:}\capt\ #2}%
      }%
   \ifdim \wd\@tempboxa >\captionwidth {%
      \rightskip=\@captionmargin\leftskip=\@captionmargin
      \unhbox\@tempboxa\par}%
   \else
      \hbox to\captionwidth{\hfil\box\@tempboxa\hfil}%
   \fi}%
\def\fignr{\small\sffamily\bfseries}%
\def\capt{\small\sffamily}%

\newdimen\@captionmargin\@captionmargin2cm\relax
\newdimen\captionwidth\captionwidth\hsize\relax

\def\eqref#1{(\protect\ref{#1})}

\def\proof{\@ifnextchar[{\@proof}{\@proof[\unskip]}}
\def\@proof[#1]{\noindent{\bf Proof #1.}\enspace}

\def\hint{\noindent Hint: }
\def\problem{\noindent{\bf Problem.} }

\def\@mt#1{\ifmmode#1\else$#1$\fi}
\def\qed{\hfill\@mt{\Box}}
\def\qqed{\hfill\@mt{\Box\enspace\Box}}

\def\cA{{\cal A}}
\def\cM{{\cal M}}
\def\cF{{\cal F}}
\def\cK{{\cal K}}
\def\cU{{\cal U}}
\def\cC{{\cal C}}
\def\cP{{\cal P}}
\def\tP{{\tilde P}}
\def\tZ{{\tilde Z}}
\def\cV{{\cal V}}
\def\fg{{\frak g}}
\def\tr{\text{tr}}
\def\cZ{{\cal Z}}
\def\cD{{\cal D}}
\def\bR{{\Bbb R}}
\def\cE{{\cal E}}
\def\bZ{{\Bbb Z}}
\def\bN{{\Bbb N}}
\def\bQ{{\Bbb Q}}

\def\bysame{\same[\kern2cm]\,}

\def\br#1{\left\lfloor#1\right\rfloor}
\def\BR#1{\left\lceil#1\right\rceil}

\def\abstractname{}

\@addtoreset {footnote}{page}

\renewcommand{\section}{%
   \@startsection
         {section}{1}{\z@}{-1.5ex \@plus -1ex \@minus -.2ex}%
               {1ex \@plus.2ex}{\large\bf}%
}
% \show\@seccntformat
\renewcommand{\@seccntformat}[1]{\csname the#1\endcsname .
\quad}

\def\bC{{\Bbb C}}
\def\bP{{\Bbb P}}

{\let\@noitemerr\relax
\vskip-2.7em\kern0pt\begin{abstract}
\noindent{\bf Abstract.}\enspace
We prove a finiteness property of the values of the skein polynomial 
of homogeneous knots which allows to establish large classes of
such knots to have arbitrarily unsharp Bennequin inequality (for
the Thurston-Bennequin invariant of any of their Legendrian embeddings
in the standard contact structure of $\bR^3$), and a give a short proof
% (for this class of knots) of a result of Birman and Menasco
that there are only finitely many among these knots that have given
genus and given braid index.\\[1mm]
\end{abstract}
}

\section{Introduction}

In this paper, we will show the following result on the skein
(or HOMFLY) polynomial \cite{HOMFLY}.

\begin{theorem}\label{tm}
The set
\[
\{\,P_K\,:\,\spn_l P_K\le b,\,\tl g(K)\le g\,\}
\]
is finite for any natural numbers $g$ and $b$.
\end{theorem}

Here $\spn_l P_K$ if the span of the HOMFLY polynomial
$P_K=P(K)$ of a knot $K$ in the (non-Alexander) variable $l$,
that is, the difference between its minimal and maximal degree in
$l$, $\md_l(P_K)$ and $\Md_l(P_K)$. By $\tl g(K)$ we denote the
weak genus of $K$ \cite{gen1}.

The main application of this theorem is to exhibit large
%
% On the other hand the use of the HOMFLY polynomial allows to
% give an application to Legendrian knots and to exhibit large
families of knots to have arbitrarily unsharp Bennequin
inequality for any of their realizations as topological
knot type of a Legendrian knot, which simplifies and extends the
main result of Kanda \cite{Kanda} and its alternative proofs
given by Fuchs--Tabachnikov \cite{TabFuchs} and Dasbach-Mangum
\cite{DasMan}.

\begin{corr}\label{cr4}
Bennequin's inequality becomes arbitrarily unsharp on any
sequence of (Legendrian embeddings of distinct) properly obversed
(mirrored) homogeneous knots.
\end{corr}

The main tool we use for the proof of theorem \reference{tm}
is the result of \cite{gen1} and an analysis of the skein (HOMFLY)
polynomial \cite{HOMFLY}. The proof of corollary \reference{cr4}
uses the inequality, known from work of Tabachnikov
\cite{Tab,TabFuchs}, relating the Thurston--Bennequin and
Maslov (rotation) number of Legendrian knots to
the minimal degree of their skein polynomial. This inequality
suggests that one should look at knots behaving ``nicely'' with
respect to their skein polynomial. The homogeneous knots introduced
by Cromwell in \cite{Cromwell}, are, in some sense, the largest
class of such knots. These are the knots having homogeneous diagrams,
that is, diagrams containing in each connected component (block) of the
complement of their Seifert (circle) picture only crossings of the same
sign. This class contains the classes of alternating and positive/
negative knots.

The other application we give of theorem \reference{tm} is also related
to Bennequin's paper \cite{Bennequin} and the work of Birman and
Menasco building on it. In their paper \cite{BirMen}, its referee made
the observation that there are only finitely many knots of given genus
and given braid index (theorem 2). This fact came as a bi-product of
the work of the authors on braid foliations introduced in Bennequin's
paper and bases on a rather deep theory.
Here we will use our work in \cite{gen1} to give a simple,
because entirely combinatorial, proof of a generalization of this
result for homogeneous knots. In fact we show that the
lower bound for the braid index coming from the inequality of
Franks--Williams \cite{WilFr} and Morton \cite{Morton}
gets arbitrarily large for homogeneous knots of given genus
(corollary \reference{cry}). The result of \cite{BirMen} for
homogeneous knots is then a formal consequence of ours.

\begin{corr}\label{crBM}
There are only finitely many homogeneous knots $K$
of given genus $g(K)$ and given braid index $b(K)$.
\end{corr}

We should remark that the corollary will straightforwardly generalize
to links. The arguments we will give apply for links of any given
(fixed) number of components, and clearly a link of braid index $n$
has at most $n$ components.

Since this paper was originally written, more work was done
on the subject, including by Etnyre, Honda, Ng, and in particular
Plamenevskaya \cite{Plamenevskaya}. A recent survey can be found
in \cite{Etnyre}.

\section{Knot-theoretic preliminaries}

The \em{skein (HOMFLY) polynomial}\footnote{Further names I have seen
in the literature are: 2-variable Jones, Jones-Conway, LYMFHO,
FLYPMOTH, HOMFLYPT, LYMPHTOFU, \dots} $P$ is a
Laurent polynomial in two variables $l$
and $m$ of oriented knots and links and can be defined
by being $1$ on the unknot and the (skein) relation
\begin{eqn}\label{1}
l^{-1}\,P\bigl(
\diag{5mm}{1}{1}{
\picmultivecline{0.18 1 -1.0 0}{1 0}{0 1}
\picmultivecline{0.18 1 -1.0 0}{0 0}{1 1}
}
\bigr)\,+\,
l \,P\bigl(
\diag{5mm}{1}{1}{
\picmultivecline{0.18 1 -1 0}{0 0}{1 1}
\picmultivecline{0.18 1 -1 0}{1 0}{0 1}
}
\bigr)\,=\,
-m\,P\bigl(
\diag{5mm}{1}{1}{
\piccirclevecarc{1.35 0.5}{0.7}{-230 -130}
\piccirclevecarc{-0.35 0.5}{0.7}{310 50}
}
\bigr)\,.
\end{eqn}
This convention uses the variables of \cite{LickMil}, but
differs from theirs by the interchange of $l$ and $l^{-1}$.

Let $[Y]_{t^a}=[Y]_a$ be the \em{coefficient} of $t^a$ in a polynomial
$Y\in\bZ[t^{\pm1}]$. For $Y\ne 0$,
let $\cC_Y\,=\,\{\,a\in\bZ\,:\,[Y]_a\ne 0\,\}$ and define
\[
\md Y=\min\,\cC_Y\,,\quad \Md Y=\max\,\cC_Y\,,\mbox{\qquad and\qquad}
\spn Y=\Md Y-\md Y\,
\]
to be the \em{minimal} and \em{maximal degree} and \em{span} (or
\em{breadth}) of $Y$, respectively. %We will denote by $\mc_t Y=
%[Y]_{\Md Y}$ the \em{leading coefficient} of $t$ in $Y$.

Similarly one defines for $Y\in\bZ[x_1,\dots,x_n]$
the coefficient $[Y]_X$ for some monomial $X$ in the $x_i$.
For a multi-variable polynomial the coefficient may be taken
with respect only to some variables, and is a polynomial in the
remaining variables, for example $[Y]_{x_1^k}\in
\bZ[x_2,\dots,x_n]$. (Thus it must be clear as a monomial
in which variables $X$ is meant. For example, for $X=x_1^k\in\bZ[x_1]$,
the coefficient $[Y]_{X}=[Y]_{x_1^k}\in\bZ[x_2,\dots,x_n]$ is \em{not}
the same as when regarding $X=x_1^k=x_1^kx_2^0\in\bZ[x_1,x_2]$ and
taking $[Y]_{X}=[Y]_{x_1^kx_2^0}\in\bZ[x_3,\dots,x_n]$.) 
% Analogously one defines for $Y\in\bZ[x_1,\dots,x_n]$ the coefficient
% $[Y]_Q$ for some monomial $Q$ in the $x_i$, and $\Md_{x_i}Y$,
% $\mc_{x_i}Y$, etc.

We call the three diagram fragments in \eqref{1} from left to
right a \em{positive} crossing, a \em{negative} crossing and a
\em{smoothed out} crossing (in the skein sense). The smoothing
out of each crossing in a diagram $D$ leaves a collection of
disjoint circles called \em{Seifert circles}. We write $c(D)$ for the
number of crossings of $D$ and $s(D)$ for the number of its
Seifert circles.

The \em{weak genus} of $K$ \cite{gen1} is the minimal genus of all
its diagrams, and the genus $g(D)$ of a diagram $D$ we will
call the genus of the surface, obtained by applying the Seifert
algorithm to this diagram:
\[
\tl g(K)\,=\,\,\min\left\{\,g(D)=\frac{c(D)-s(D)+1}{2}\,:\,\mbox{$D$ is
a diagram of $K$}\,\right\}\,.
\]

The \em{genus} $g(K)$ of $K$ is the minimal genus of all Seifert
surfaces of $K$ (not necessarily coming from Seifert's algorithm
on diagrams of $K$). The \em{slice genus} $g_s(K)$ of $K$ is the
minimal genus of all smoothly embedded surfaces $S\subset B^4$ with
$\partial S=K\subset S^3=\partial B^4$. Clearly $g_s(K)\le g(K)\le
\tl g(K)$. $\tl g(K)$ coincides with the usual genus for many knots,
in particular knots up to 10 crossings and homogeneous knots.

The \em{braid index} $b(K)$ of $K$ is the minimal number of strings
of a braid having $K$ as its closure. See \cite{BirMen,WilFr,Morton}.

Recall, that a knot $K$ is \em{homogeneous}, if it has a
diagram $D$ containing in each connected component of the
complement (in $\bR^2$) of the Seifert circles of $D$
(called \em{block} in \cite[\S 1]{Cromwell}) only crossings of the
same sign (that is, only positive or only negative ones).
This notion was introduced in \cite{Cromwell} as a
generalization of the notion of alternating and positive knots.

\section{On the genus and braid index of homogeneous knots}

The main tool we use for the proof of theorem \reference{tm}
is the result of \cite{gen1}. 

\begin{theorem}\label{thm}(\cite{gen1}) 
Knot diagrams of given genus (with no nugatory crossings and modulo
crossing changes) decompose into finitely many equivalence classes
under flypes \cite{MenThis} and (reversed) applications
of antiparallel twists at a crossing
\begin{eqn}\label{move}
\diag{7mm}{1}{1}{
    \picmultivecline{0.12 1 -1.0 0}{0 0}{1 1}
    \picmultivecline{0.12 1 -1.0 0}{1 0}{0 1}
}\llra\quad
\diag{7mm}{3}{2}{
  \picPSgraphics{0 setlinecap}
  \pictranslate{0.5 1}{
    \picrotate{-90}{
      \lbraid{0 -0.5}{1 1}
      \lbraid{0 0.5}{1 1}
      \lbraid{0 1.5}{1 1}
      \pictranslate{-0.5 0}{
      \picvecline{0.03 1.95}{0 2}
      \picvecline{0.03 -.95}{0 -1}
    }
    }
    }
}
\,.
\end{eqn}
\end{theorem}

This theorem allows to define for every natural number $g$ an integer
$d_g$ as follows (see \cite{gen1} for more details): call 2 crossings
in a knot diagram equivalent, if there is a sequence of flypes making
them to form a clasp
$
\diag{6mm}{2}{1}{
  \picrotate{-90}{
    \lbraid{-0.5 0.5}{1 1}
    \lbraid{-0.5 1.5}{1 1}
    \picvecline{-0.95 1.9}{-1 2}
    \picvecline{-0.05 0.1}{0 0}
  }
}\,,
$ in which the strands are reversely oriented.
One checks that this is an equivalence relation. Then $d_g$ is the
maximal number of equivalence classes of crossings of diagrams of genus
$g$. The theorem ensures that $d_g$ is finite. It follows from the work
of Menasco and Thistlethwaite \cite{MenThis} that $d_g$ can we
expressed more self-containedly as
\[
d_g\,=\,1+\sup\left\{\,i\in\bR\,:\,\limsup_{n\to\infty}
\,\frac{a_{n,g}}{n^i}\,>\,0\,\right\}\,,
\]
where $a_{n,g}$ is the number of alternating knots of $n$ crossings
and genus $g$. (Note that it is not \em{a priori} clear, whether the
supremum on the right is integral or even finite.)

\proof[of theorem \reference{tm}] It follows from theorem
\reference{thm} that we can w.l.o.g. consider only one equivalence
class $\cD$ of diagrams of genus $g$ modulo the move \eqref{move}.
We will now argue that
the skein relation for the HOMFLY polynomial implies that
for a knot diagram $D$ in $\cD$ we have for its polynomial $P(D)=P_D$ 
\begin{eqn}\label{ff}
(l^2+1)^{d_g}P(D)\,=\,\sum_{i=1}^{n_\cD}\,l^{t_i(D)}\, L_{i,\cD}\,,
\end{eqn}
where the number $n_\cD$ and the polynomials $L_{i,\cD}\in\bZ[l^{\pm 1},
m^{\pm 1}]$ depend on $\cD$ only, and the only numbers depending on
$D$ are the $t_i$. (The reader may compare to
\cite[proof of theorem 3.1]{gwg}
for the case of the Jones polynomial, which is analogous.)
An easy consequence of this skein relation \eqref{1} is the relation 
\begin{eqn}\label{g}
P_{2n+1}\,=\,\frac{(il)^{2n}-1}{l+l^{-1}}m\,P_\infty+(il)^{2n}P_1\,,
\end{eqn}
where $L_{2n+1}$ is the the link diagram obtained by $n$ $\bt$ moves
at a positive crossing $p$ in the link diagram $L_1$ (we call the
tangle in $L_{2n+1}$ containing the $2n+1$ crossings so obtained a
``twist box''), $L_\infty$ is $L_1$ with $p$ smoothed out, and $P_i$ is
the polynomial of $L_i$. (Compare also to the formula (8) of \cite{gwg},
but note that this formula has a misprint: the second term on
the right must be multiplied by $t$.) To obtain \eqref{ff},
iterate \eqref{g} over the at most $d_g$ twist boxes of a genus $g$
diagram. The factor $(l^2+1)^{d_g}$ is used to get disposed of the
denominators. 

{}From \eqref{ff} we obtain for a diagram $D$ of genus $g$
\begin{eqn}\label{eq1}
\big|\,\bigl[\,(l^2+1)^{d_g}P(D)\,\bigr]_{l^pm^q}\,\big|\,\le\,C_{g}
\end{eqn}
for some constant $C_g$ depending on $g$ only.

Morton showed in \cite{Morton} that $[P(D)]_{l^pm^q}=0$ if $q>2g(D)$,
and, as well-known, the same is true for $q<0$ (we assume that $D$ is
a \em{knot} diagram). Furthermore, 
it follows from the identity $P_D(l,-l-l^{-1})=1$ that 
$[P(D)]_{l^pm^q}\ne 0$ for some $|p|\le 2g(D)$. Thus, if
$\spn_l P_D\le b$ and $g(D)\le g$, then $[P_D]_{l^pm^q}=0$
for $|p|+q>C'_{b,g}$ with $C'_{b,g}$ depending only on $b$ and $g$, and
hence the same is true for $(l^2+1)^{d_g}P_D$. But we already saw that
$(l^2+1)^{d_g}P_D$ has uniformly bounded coefficients \eqref{eq1},
so that
\[
\{\,(l^2+1)^{d_g}P_D\,:\,\spn_l P_D\le b,\,g(D)\le g\,\}
\]
is finite. From this the theorem follows because multiplication
with $(l^2+1)^{d_g}$ is injective (the polynomial ring is an
integrality domain). \qed

\begin{corr}\label{cry}
There are only finitely many homogeneous knots $K$
of given genus $g(K)$ and given value of $\spn_l P_K$.
\end{corr}

We should point out that (trivially) a given knot may have infinitely
many diagrams of given genus, and that even infinitely many different
knots may have diagrams of given genus with the same HOMFLY polynomial
\cite{Kanenobu}. This, not unexpectedly, shows that the combinatorial
approach has its limits.

\proof[of corollary \reference{cry}] Combine theorem \reference{tm}
with the facts that for a homogeneous knot $K$ we have $g(K)=\tl g(K)$
\cite[corollary 4.1] {Cromwell}, that there are only finitely many
homogeneous knots already of given Alexander polynomial \cite[corollary
3.5]{gen1}. \qed

Finally, as the inequality $\max\deg_mP(K)\le 2\tl g(K)$ of
Morton \cite{Morton} is known to be sharp in very many cases, we
are led to conjecture more.

\begin{conjecture}
The set
\[
\{\,P_K\,:\,\mbox{$K$ knot, } \spn_l P_K\le b,\,\max\deg_mP_K\le g\,\}
\]
is finite for any natural numbers $g$ and $b$.
\end{conjecture}

\section{The HOMFLY polynomial and Bennequin's inequality for
Legendrian knots}

A \em{contact structure} on a smooth 3-manifold is a 1-form $\al$ with
$\al\wedge d\al\ne 0$ (which is equivalent to the non-integrability 
of the plane distribution defined by $\ker\al$). In the following we
consider the $1$-form $\al=dx+y\,dz$ on $\bR^3(x,y,z)$, called
the \em{standard contact space}. A \em{Legendrian
knot} is a smooth embedding $\cK:S^1\to \bR^3$ with $\al\left(\frac{%
\partial \cK}{\partial t}\right)\equiv 0$. Each such knot has its
underlying topological knot type $K=[\cK]$ and two fundamental
invariants in contact geometry known as the \em{Thurston--Bennequin
number}
$tb(\cK)$ and \em{Maslov index} $\mu(\cK)$. (See \cite{TabFuchs,Ferrand}
for an excellent introductory account on this subject.)

\begin{defi}
The Thurston--Bennequin number $tb(\cK)$ of a Legendrian knot $\cK$
in the standard contact space 
is the linking number of $\cK$ with $\cK'$, where $\cK'$ is obtained
from $\cK$ by a push-forward along a vector field transverse to the
(hyperplanes of the) contact structure.

The Maslov (rotation) index $\mu(\cK)$ of $\cK$ is the degree of the map
\[
t\in S^1\,\mapsto\,\frac{\pr\,\frac{\partial \cK}{\partial t}(t)}{
\bigl |\pr\,\frac{\partial \cK}{\partial t}(t)\bigr |}\in S^1\,,
\]
where $\pr\,:\,\bR^3\to\bR^2\simeq \bC$ is the projection
$(x,y,z)\mapsto(y,z)$.
\end{defi}

Both invariants $tb(\cK)$ and $\mu(\cK)$ can be interpreted
in terms of a regular diagram of the (topological) knot $[\cK]$, and 
thus it was recently realized that the theory of polynomial invariants
of knots and links in $\bR^3$, developed after Jones \cite{Jones},
can be applied in the context of Legendrian knots to give inequalities
for $tb$ and $\mu$. In particular we have the inequality
\begin{eqn}\label{tbm} 
tb(\cK)+|\mu(\cK)|\,\le\,\md_lP([\cK])-1\,.
\end{eqn}
This follows from the work of Morton \cite{Morton} and Franks--Williams
\cite{WilFr}, and was translated to the Legendrian knot context
by Tabachnikov and Fuchs \cite{TabFuchs}. See also
\cite{Tab,CGM,GorHill,Ferrand}.

On the other hand, a purely topological inequality was previously known
for a while~-- Bennequin's inequality. In \cite{Bennequin}, Bennequin
proved
\begin{eqn}\label{bi} 
tb(\cK)+|\mu(\cK)|\,\le\,2g([\cK])-1\,.
\end{eqn}
This inequality was later improved by Rudolph \cite{Rudolph3}
who showed
\begin{eqn}\label{rbi} 
tb(\cK)+|\mu(\cK)|\,\le\,2g_s([\cK])-1\,,
\end{eqn}
where $g_s(K)$ is the slice (4-ball) genus of $K$. This improvement
used the proof of the Thom conjecture by Kronheimer and Mrowka,
achieved originally by gauge theory \cite{KroMro,KroMro2},
and later much more elegantly by Seiberg--Witten invariants
\cite{KroMro3}. 

While the r.h.s. of \eqref{bi} and \eqref{rbi} are invariant w.r.t.
taking the mirror image, the l.h.s. are strongly sensitive, so we have
\begin{eqn}
\tau'(L)\le \tau(L)\le 2g_s(L)-1\le 2g(L)-1
\end{eqn}
for any topological knot type $L$, where
\[
\tau'(L):=\max\{tb(\cK)+|\mu(\cK)|\,:\,[\cK]=L\,\}
\]
and
\[
\tau(L):=\max\{tb(\cK)+|\mu(\cK)|\,:\,[\cK]\in\{L,!L\}\,\}\,=\,
\max(\tau'(L),\tau'(!L))\,,
\]
and $!L$ is the obverse (mirror image) of $L$.

In \cite{Kanda}, Kanda used an original argument and the theory of
convex surfaces in contact manifolds developed mainly by Giroux
\cite{Giroux} to show that the inequality $\tau'\le 2g-1$ can get
arbitrarily unsharp, i.e.
$\exists\{L_i\}\,:\,\tau'(L_i)-2g(L_i)\to-\infty$.
(Here, and in the following, an expression of the form `$x_n\to\infty$'
should abbreviate $\lim\limits_{n\to\infty}x_n=\infty$. Analogously
`$x_{n_m}\to\infty$' should mean the limit for $m\to\infty$ etc.)
In Kanda's paper, all $L_i$ are alternating pretzel knots, and hence of 
genus 1, so that for these examples in fact we also have
$\tau'(L_i)-2g_s(L_i)\to-\infty$.

It was realized (see the remarks on \cite[p. 1035]{TabFuchs}) that
Kanda's result admits an alternative proof using \eqref{tbm} (whose
proof in turn is also ``elementary'' in a sense discussed more
detailedly in \cite{Ferrand}). Other examples (connected sums of
two $(2,\,.\,)$-torus knots) were given by Dasbach and Mangum \cite
[\S 4.3]{DasMan}, for which even $\tau-2g\to-\infty$. However, their
examples do not apply for the slice version \eqref{rbi} of
Bennequin's inequality. In \cite{Ferrand} it was observed that
Kanda's result also follows from the work of Rudolph \cite{Rudolph,%
Rudolph2}. 

Here we give a larger series of examples of knots with $2g-\tau
\to\infty$ containing as very special cases the previous ones
given by Kanda and Dasbach--Mangum. These knots show that
the inexactness of Bennequin's inequality is
by far not an exceptional phenomenon. While arguments also
use \eqref{tbm} (and hence are much simpler than the original
proof of Kanda), they still apply in many cases also for the
slice version \eqref{rbi} of Bennequin's inequality. Similar
reasoning works for links of any fixed number of components, but
for simplicity we content ourselves only with knots.

From theorem \reference{tm}, the aforementioned application to
the unsharpness of Bennequin's inequality is almost straightforward.
We formulate the consequence somewhat more generally and more
precisely than in the introduction.

\begin{theorem}\label{cr3}
Let $\{L_i\}$ be a sequence of knots, such
that only finitely many of the $L_i$ have the same skein polynomial.
Then $2\tl g(L_i)-\min(\tau'(L_i),\tau'(!L_i))\to\infty$.
If additionally $\tl g(L_i)\le C$ for some constant $C$, then even
$\min(\tau'(L_i),\tau'(!L_i))\to-\infty$.
\end{theorem}

The condition $g=\tl g$ is very often satisfied, but unfortunately this
is not always the case, as pointed out by Morton \cite[remark 2]
{Morton}. Worse yet, as shown in \cite{gwg}, there cannot be
any inequality of the type $\tl g(K)\le f(g(K))$ for any
function $f\,:\,\bN\to\bN$ for a general knot $K$. Nevertheless,
by the results mentioned in the proof of corollary \ref{cry},
any sequence of homogeneous knots satisfies $g(L_i)=\tl g(L_i)$
and the condition of theorem \reference{cr3}. In particular, we have

\begin{corr}\label{cr5}
If $\{L_i\}$ are negative or achiral homogeneous knots,
then $2g(L_i)-\tau'(L_i)\to
\infty$. If $L_i$ are negative or additionally $g(L_i)\le C$ for some 
constant $C$, then even $2g_s(L_i)-\tau'(L_i)\to\infty$. \qed
\end{corr}

\begin{rem}%(
\def\labelenumi{\theenumi)}
Before we prove theorem \reference{cr3}, we make some comments
on corollary \reference{cr5}.
\begin{enumerate}
\item Clearly for an achiral knot $L$ we have $\tau(L)=\tau'(L)$,
so that in the case all $L_i$ are achiral (like the examples $T_{2,n}
\#T_{2,-n}$, with $T_{2,n}$ being the $(2,n)$-torus knot, given in
\cite{DasMan}) the stronger growth statement with $\tau'$
replaced by $\tau$ holds, $2g-\tau\to\infty$.
\item Contrarily, the statement $2g-\tau\to\infty$ is not true in the
negative case: Tanaka \cite[theorem 2]{Tanaka} showed that $\tau'=2g-1$
for positive knots.
% \footnote{These are the mirror images of positive
% knots. Note, that Tanaka uses a contact structure which differs from
% ours by
% a reflection, so we need to use the obverses of the knots he
% considers.}.
On the other hand, this means that for negative knots
$2g-\tau'\to\infty$, and in fact $2g_s-\tau'\to\infty$, as by \cite{pos}
$g=g_s$ for positive (and hence also for negative) knots. However,
we have from \cite{beha} the
stronger statement that $\tau'\to-\infty$, which also holds for
almost negative knots (see \cite[\S 5]{apos}).
\item The conditions can be further weakened. For example we
can replace achirality by self-conjugacy of the HOMFLY polynomial
(invariance under the interchange $l\Lra l^{-1}$) and
`negative' by `$k$-almost negative' for any fixed number $k$,
as the condition $g=\tl g$ in corollary \reference{cr3} can in fact be
weakened to $\tl g\le f(g)$ for any (fixed) function $f\,:\,\bN\to\bN$.
However, in latter case 
the assumption needs to be retained that only finitely many $L_i$
have the same polynomial. (This is known to be automatically
true for $k\le 1$ \cite{beha,apos}, but not known for $k\ge 2$.)
\item The fact that our collection of examples is richer than the one
of Kanda can be made precise followingly: the number of all pretzel
knots of at most $n$ crossings is $O(n^3)$, while it follows from
\cite{ErnSum} and \cite{fib} that the number of achiral and
positive knots of crossing number at most $n$, already 
among the 2-bridged ones, grows exponentially in $n$.
%
% \item It was known, originally from \cite{CochGom}, that the
% properties negative and achiral are mutually exclusive. In
% fact, achirality can be regarded to be ``in the middle'' between
% positivity and ngeativity. From our proof it will get clear that
% the argument should work for many other knots which are ``inbetween''
% negative and achiral, but a self-contained definition of what this
% ``between'' should mean seems lacking. It has to be related to the
% (partial) order of knots defined by switching positive crossings to
% negative (which in turn is related to the one of Cochran and Gompf
% \cite[proposition 2.2]{CochGom}), but in some non-straightforward way.
%
\item The boundedness condition on the genus in the achiral case
is essential (at least for this method of proof) as show the examples
$T_{2,n}\#T_{2,-m}$ of Dasbach and Mangum, on which the skein
polynomial argument fails for the slice genus.
\end{enumerate}
\end{rem}

\proof[of theorem \reference{cr3}] If $\tl g(L_i)\le C$, then theorem
\reference{tm} implies that $\Md_lP(L_i)-\md_lP(L_i)=\spn_l P
(L_i)\to\infty$. Thus $\min(\md_lP(L_i),\md_lP(!L_i))\to-\infty$, and
the assertions follow from \eqref{tbm} and \eqref{rbi}. Else there is a
subsequence $\{L_{i_j}\}$ with $\tl g(L_{i_j})\to\infty$. Clearly,
\[
\min(\md_lP(L_{i_j}),\md_lP(!L_{i_j}))\le 0\,,
\]
so that $2\tl g(L_{i_j})-
\min(\tau'(L_{i_j}),\tau'(!L_{i_j}))\to\infty$, that is, $\{L_i\}$
always has a subsequence with the asserted property. Applying the
argument on any subsequence of $\{L_i\}$ gives the property on the
whole $\{L_i\}$. \qed

As a final remark, there is another inequality, proved in
\cite{ChmGor} and \cite{Tab}, involving the
Kauffman polynomial $F$ (in the convention %\em{conjugate} to that
of \cite{Kauffman}),
\begin{eqn}
tb(\cK)\,\le\,-\Md_aF([\cK])-1\,.
\end{eqn}
It gives in general better estimates on $tb(\cK)$, but lacks the
additional term $|\mu(\cK)|$ and also a translation to the transverse
knot context (see \cite[remark at end of \S 6]{Ferrand}). Contrarily,
\eqref{tbm}
admits a version for transverse knots as well (in which case the term
$|\mu(\cK)|$ is dropped; see \cite{GorHill} and \cite[theorem 2.4]{TabFuchs}),
and so our results hold in the transverse case, too. Moreover, as
remarked by Ferrand in \cite[\S 8]{Ferrand} (see also \cite[Problem,
p.\ 3428]{Tanaka}), the inequality
\begin{eqn}\label{PF}
-\Md_aF(K)\le \md_lP(K)
\end{eqn}
is not always satisfied. Among the 313,230 prime knots of at most 15
crossings tabulated in \cite{KnotScape} there are 134 knots $K$ such
that at least one of $K$ and $!K$ fails to satisfy \eqref{PF}. The
simplest examples are two 12 crossing knots, one of them, $12_{1584}$,
being quoted by Ferrand.

\noindent{\bf Acknowledgement.} I wish to thank to E. Ferrand and
S. Tabachnikov for some helpful conversations and for introducing
me to the subject. I also wish to thank to K. Hulek for organizing
and inviting me to the conference ``Perspectives in Mathematics''
and letting me there prepare a part of this manuscript.

{\small
\let\old@bibitem\bibitem
\def\bibitem[#1]{\old@bibitem}

}

\end{document}